\documentclass[11pt]{amsart}
\usepackage{geometry}                % See geometry.pdf to learn the layout options. There are lots.
\usepackage{graphicx}
\usepackage{amssymb}
\usepackage{epstopdf}
\usepackage{color}
\title{Appendix to Ólimits for partial maxima of gaussian random vectorsÓ }

\author{James Kuelbs}
\address{James Kuelbs\\ Department of Mathematics,  University of Wisconsin, Madison, WI 53706-1388} 
\email{kuelbs@math.wisc.edu}

\author{Joel Zinn}
\address{Joel Zinn\\ Department of Mathematics, Texas A\&M University, College Station, TX 77843-3368} 
\email{joel.zinn@gmail.com}

\begin{document}
\maketitle
%\section{}
%\subsection{}
{\bf Abstract.} This appendix provides a short proof for sample path continuity of the Brownian motion induced by an arbitrary centered Gaussian measure on a separable Banach space, and also some perturbation results for the spectrum of compact self-adjoint operators on a Hilbert space.
 \bigskip

{\bf Introduction.} Let $\gamma$ be a non-degenerate mean zero Gaussian measure on the separable Banach space $E$ with norm $q$. Here we establish sample function continuity of the $\gamma$-generated Brownian motion used in [3]. This result appeared earlier in [2], but the proof we provide is direct and fairly short so perhaps it could be of independent interest. 

A second topic in this appendix deals with results for compact self-adjoint operators on the Hilbert space $H$. Here an important reference is [5], and although we are unaware of the results presented in Theorem 3 below,  much of what is established already follows from results appearing in [5].

{\bf 1. Sample Function Continuity.} As above $\gamma$ is a mean zero Gaussian measure on $E$ with norm $q$. In addition, assume $\Omega_E$ is the space of continuous functions $x$ from $[0,\infty)$ into $E$ such that $x(0)=0$, and $\mathcal{F}$ is the $\sigma$-field of $\Omega_E$ generated by the functions $x \rightarrow x(t), 0 \le t < \infty$. Our next lemma provides a proof that there exists a probability measure $P_{\gamma}$ on $(\Omega_E,\mathcal{F})$ such that if $0=t_0<t_1\cdots< t_n$, then the random vectors
$$
x(t_j)-  x(t_{j-1}), j=1,\cdots ,n, ~{\rm{are~independent}}, \leqno(1.1)
$$
and $x(t_j)- x(t_{j-1})$ has distribution $\gamma_{t_j - t_{j-1}}$ on $E$, where $\gamma_{s}(A)= \gamma(A/\sqrt{s})$ for Borel subsets $A$ of $E$ when $s>0$ and $\gamma_0=\delta_0.$ In particular, the stochastic process  $\{W(t): t\ge 0\}$ defined on $(\Omega_E, \mathcal{F},P_{\gamma})$ by $W(t,x)=x(t)$ has stationary independent mean zero Gaussian increments. We will call it $\gamma$-Brownian motion on $E$. 

{\bf Proposition 1.}
Let $\gamma$ be a mean zero Gaussian measure on $E$. 
Then, the $E$-valued Brownian motion $W=\{W(t): t\ge 0\}$ defined on $(\Omega_E, \mathcal{F},P_{\gamma})$ by $W(t,x)= x(t)$ exists. In particular, it is sample path continuous and mean zero with all finite dimensional distributions jointly Gaussian with\\ 
\indent {(i) $W(0)=0,$ }\\
\indent {(ii) $W$ has stationary independent mean zero Gaussian increments as indicated above, and}\\
\indent {(iii) if the support of $\gamma$} is a closed subspace $F$ of $E$, then the probability $P_{\gamma}$ has support on $\Omega_F$, the space of continuous functions on $[0,\infty)$ with values in $F$, and $W$ is a $\gamma$ Brownian motion on F.\\

{\bf Proof.} For sample function continuity it suffices to show that $\{W(t): t \ge 0\}$ is sample continuous on $[0,b]$ for every integer $b \ge 1$.

Let $G, G_1,G_2,\cdots$ be i.i.d. $E$-valued random vectors with distribution $\gamma$ defined on the probability space $(\Omega, \mathcal{G}, Q)$, and set $S_0=0$ and $S_k= \sum_{j=1}^k G_j$ for $k \ge 1$. Then, for the integer $b\ge 1$ fixed and every integer $N\ge 1$ we define on $(\Omega,\mathcal{G},Q)$ the piecewise linear sample continuous processes
$$
X_N(t)(\omega)=S_k(\omega)/\sqrt{2^N}, t= k/2^N, k=0,1,\cdots,b2^N,
$$
and $X_N(t)(\omega)$ is linear and continuous on each of the subintervals $I_k= [t_{k-1}, t_k], 1 \le k \le b2^N$ where $t_0=0$ and $t_k=k/2^N, 1 \le k \le b2^N$. Also define for $x$ continuous and $E$-valued on $[0,b]$ the maps 
$$
\Pi_N(x)(t)=x(k/2^N),t=k/2^N,  k=0,\cdots,  b2^N, N\ge 1,
$$
where $\Pi_N(x)(\cdot)$ is piecewise linear and continuous on the subintervals $I_k=[t_{k-1},t_k]$ with $t_0=0, t_k= k/2^N, 1\le k \le b2^N.$ Then the maps $\Pi_N$ take the separable Banach space $C_E[0,b]$ of $E$-valued continuous functions $x$ with norm $\sup_{t \in [0,b]}q(x(t))$ into itself, and for each Borel subset $A$ of $C_E[0,b]$ we define for $N\ge 1$ the probability measures
$$
P_N(A)= Q( \{\omega: X_N(\cdot)(\omega) \in A\}).
$$
Then, it is easy to see that $P_N$ has support on $\Pi_N(C_E[0,b])$, 
and for $N\ge 2$
$$
P_{N-1}= P_N^{\Pi_{N-1}}.
$$
Letting $\rho(\cdot, \cdot)$ denote the Prokhorov metric on the Borel probability measures of $C_E[0,b]$, we then have from Lemma 1.2 of [4]
$$
\rho(P_{N-1},P_N)= \rho(P_N^{\Pi_{N-1}},P_N)\le \max( \epsilon,\delta), \leqno(1.2)
$$
provided 
$$
P_N(x \in C_E[0,b]: \sup_{ t \in [0,b]}q(x(t) -\Pi_{N-1}x(t)) \ge \delta)<\epsilon.
$$
Now
$$
P_N(x \in C_E[0,b]: \sup_{ t \in [0,b]}q(x(t) -\Pi_{N-1}x(t)) \ge \delta)=Q(\sup_{t \in [0,b]}q(X_N(t) - \Pi_{N-1}(X_N(t))) \ge \delta),
$$
and since
$$
Q(\sup_{t \in [0,b]}q(X_N(t) - \Pi_{N-1}(X_N(t))) \ge \delta) {\color{blue} =} Q( \sup_{t \in I_k,1\le k \le b2^N}q(X_N(t) - \Pi_{N-1}(X_N(t))) \ge \delta),
$$
with
$$
 \sup_{t \in I_k,1\le k \le b2^N}q(X_N(t) - \Pi_{N-1}(X_N(t)))
$$
$$ 
\le \sup_{1\le k \le b2^N}2q(X_N(\frac{k}{2^N})- X_N(\frac{k-1}{2^N}))=\sup_{1 \le k\le b2^N}2q(G_k/\sqrt{2^N}),
$$
we have
$$
Q(\sup_{t \in [0,b]}q(X_N(t) - \Pi_{N-1}(X_N(t))) \ge \delta) \le b2^N Q(q(G)\ge \sqrt{2^N}\delta/2)
$$
$$
\le b2^{N(1-\frac{r}{2})}(2/\delta)^{r} \int_{\Omega}q^{r}(G)dQ.
$$
Since $q(G)$ has finite $r^{th}$-moment for all $r>0$ , for simplicity take $r=4$ and $\delta= 2^{-N/8}$, which implies 
$$
Q(\sup_{t \in [0,b]}q(X_N(t) - \Pi_{N-1}(X_N(t))) \ge 2^{-N/8})\le 16b \int_{\Omega}q^4(G)dQ 2^{N(1-2+\frac{1}{2})},
$$
and hence from (1.2) that for all $N$ sufficiently large
$$
\rho(P^N,P^{N-1}) \le 2^{-N/8},
$$
Therefore, $\{P^N: N\ge 1\}$ is a Cauchy sequence  of probability measures on the Borel subsets of $C_E[0,b]$ and since the space of all Borel probability measures on $C_E[0,b]$ is complete in the Prokhorov metric, see [4, Theorem 1.11], the sequence $\{P^N\}$ converges weakly to a Borel probability measure $P$ on $C_E[0,b]$. Observe that the probability measures $P^N$ are such that for $t \in [0,b]$, a dyadic rational,
$$
\mathcal{L}(X_N(t)) = \gamma_{t}  
$$
for all $N$ sufficiently large. Similarly, disjoint intervals with dyadic rational endpoints are eventually independent for $P^N$ for all $N$ sufficiently large and hence the increments over disjoint intervals with dyadic rational endpoints are also independent for $P$ with laws as indicated in the definition of $\gamma$-Brownian motion.  Finally, using the fact that the sample functions are continuous, the increments over arbitrary disjoint intervals are limits of increments over intervals with dyadic end points, so it follows easily using characteristic functions that the increments are independent over arbitrary disjoint intervals of $[0,b]$ with respect to $P$ and have laws as required. That the finite dimensional distributions are mean zero and Gaussian easily follows from the fact that the process has independent mean zero Gaussian increments, and (iii) is obvious from the construction used to prove the sample function continuity. Since the integer $b\ge 1$ was  arbitrary, the proposition is proved.
\bigskip

{\bf 2.1. Eigenvalue Comparison for Self-Adjoint Compact Operators.} Let $A$ be a bounded linear operator on a Hilbert space $H$. The complex number  $\mu$ is an eigenvalue of $A$ if there is an $x \in H, x \not=0,$ such that 
$$
Ax=\mu x.
$$
Each such non-zero $x$ is said to be an eigenvector of $A$ corresponding to the eigenvalue $\mu$. Note that $\mu=0$ is possibly an eigenvalue of $A$, and we asume $H$ is infinite dimensional throughout. If $H$ is finite dimensional the results in [1], pp 141-162, establish an analogue of what we describe below.

Let $A$ be a compact, self-adjoint operator on a Hilbert space $H$ with $||A|| \not= 0$. Then, $A$ has at least one non-zero eigenvalue, and all eigenvalues of $A$ are real numbers. Moreover, $A$ has a finite or denumerable number of non-zero eigenvalues $\{\mu_k\},$ and $\mu_k \rightarrow 0$ provided there are infinitely many non-zero eigenvalues. Eigenvectors corresponding to distinct eigenvalues are orthonormal, and each $\mu_{k} \not= 0$ has at most finitely many eigenvectors (orthogonal to each other) and is repeated as many times as there are distinct corresponding orthogonal eigenvectors. Every element of the form $Ax$ can be written as
$$
Ax= \sum_{k \ge 1} (Ax,\phi_k)\phi_k = \sum_{k \ge 1}\mu_k(x,\phi_k)\phi_k, \leqno(2.1.1)
$$
where $\phi_k$ are the orthogonal eigenvectors corresponding to non-zero eigenvalues. In order that the orthogonal eigenvectors $\{\phi_k\}$ corresponding to non-zero eigenvalues form a complete orthonormal system for $H$ it is necessary and sufficient the $A$ does not have zero as an eigenvalue, that is, $Ax \not= 0$ when $x \not= 0$. Of course, it is always possible to take the eigenvectors corresponding to the eigenvalue zero to be an orthonormal basis of the null space of $A$, so that when these eigenvectors are combined with those corresponding to the non-zero eigenvalues $\{\mu_k\}$, the union forms an orthonormal basis of $H$.

To determine the eigenvalue (non-zero) with largest  absolute value and a corresponding eigenvector we observe that the compactness of $A$ implies the extremal problem 
$$
|(Af,f)|= \sup_{||x||=1} |(Ax,x)|
$$
admits solutions of norm one, all such solutions $f=\phi$ are eigenvectors of $A$ (see pages 231-2 of [5]), and the corresponding eigenvalue $\mu_1$  of $\phi_1$ is such that
$$
|\mu_1|= \sup_{||x||=1} |(Ax,x)|.
$$
Also, note from page 230 of [5] that
$$
||A||= \sup_{||x||=1}|(Ax,x)|
$$
and hence if the supremum is zero, this implies $A=0$ and there are no non-zero eigenvalues. Of course, there are then many orthogonal eigenvectors corresponding to the eigenvalue zero.

To find a second eigenvector for  $A$ corresponding to a non-zero eigenvalue solve
$$
|(Af,f)|= \sup_{||x||=1, (x,\phi_1)=0} |(Ax,x)|,
$$
assuming
$$
 \sup_{||x||=1, (x,\phi_1)=0} |(Ax,x)|\not= 0.
$$
If this supremum is zero, $A=0$ on \{$ x \in H: (x,\phi_1)=0\}$, there are no additional non-zero eigenvalues, and any orthonormal basis of $\{x \in H: (x,\phi_1)=0\}$ consists of orthonormal eigenvectors for the eigenvalue zero.
When the supremum is non-zero, compactness of $A$ again implies there is an orthonormal eigenvector $\phi_2$ orthogonal to $\phi_1$,  and the corresponding eigenvalue $\mu_2$ is such that
$$
|\mu_2|= \sup_{||x||=1, (x,\phi_1)} |(Ax,x)|.
$$
Once we have found the orthogonal eigenvectors $\phi_,\cdots,\phi_{n-1}$ we solve
$$
|(Af,f)|= \sup_{||x||=1, (x,\phi_j)=0, j=1,\cdots, n-1} |(Ax,x)|.
$$
Again compactness of $A$ implies there is an orthonormal eigenvector $\phi_n$ orthogonal to $\phi_1,\cdots,\phi_{n-1}$,  and the corresponding eigenvalue $\mu_n$ is such that
$$
|\mu_n|= \sup_{||x||=1, (x,\phi_j), j=1.\cdots,n-1} |(Ax,x)|.
$$

This process stops after a finite number of iterations if some $\mu_k=0$, or it proceeds indefinitely obtaining a sequence of orthonormal eigenvectors with nonzero eigenvalues of $A$, each appearing as many times as its multiplicity indicates. As mentioned previously, in case there are infinitely many  eigenvalues $\{\mu_k\}$, then the compactness of $A$ forces $\mu_k$ to converge to zero so at most there are countably many eigenvectors that correspond to nonzero eigenvalues, and (2.1.1)  
holds for all $x \in H$.

To this point the eigenvalues $\{ \mu_k\}$ have been ordered in terms of their absolute values with $|\mu_k|$ non-increasing. Now we arrange the (strictly) positive and (strictly) negative terms of the sequence  $\{ \mu_k\}$ into separate sequences
$$
\mu_1^{+}, \mu_2^{+}, \ldots ~{\rm{and}}~~ \mu_1^{-}, \mu_2^{-}, \ldots, \leqno(2.1.2)
$$
where the $\mu_k^{+}$ are non-increasing and the $\mu_k^{-}$ are non-decreasing, and each appears as many times as the multiplicity of its corresponding eigenvectors dictates. Each of these sequences can be finite, infinite, or even empty, and we denote the corresponding orthonormal eigenvectors by
$$
\phi_1^{+}, \phi_2^{+}, \ldots ~{\rm{and}}~~ \phi_1^{-}, \phi_2^{-}, \ldots. \leqno(2.1.3)
$$
From (2.1.1)
we then have for each $x\in H$ that
$$
(Ax,x)= \sum_{k\ge 1}\mu_k^{+}(x,\phi_k^{+})^2~ + ~  \sum_{k\ge 1}\mu_k^{-}(x,\phi_k^{-})^2, \leqno(2.1.4)
$$
and it is immediate to see
$$
\mu_n^{+}= \sup_{||x||=1, (x,\phi_k^{+})=0, k=1,\cdots,n-1} (Ax,x)= (A\phi_n^{+}, \phi_n^{+}) \leqno(2.1.5)
$$
and
$$
\mu_n^{-}= \inf_{||x||=1, (x,\phi_k^{-})=0, k=1,\cdots,n-1} (Ax,x)= (A\phi_n^{-}, \phi_n^{-}). \leqno(2.1.6)
$$

This characterization of the $n^{th}$ eigenvalue and the $n^{th}$ eigenvector depends on knowing the previous $n-1$ quantities. The next result is is a direct characterization which was obtained for finite dimension matrices (quadratic forms) by E. Fisher in 1905 and extended by D. Hilbert in 1912. These references, as well as much of these notes through (2.1.14)  
below, can be found in [5], pages 231-240.
\bigskip

{\bf Theorem 1.} Let $h_1,\ldots,h_{n-1}$ be  $n-1$ elements of $H$, and
$$
\nu(h_1,\ldots,h_{n-1})= \sup_{||x||=1, (x,h_j)=0, j=1,\ldots,n-1} (Ax,x). \leqno(2.1.7)
$$
If there are $n$  strictly positive eigenvalues, then
$$
\inf_{h_1,\ldots,h_{n-1}} \nu(h_1,\ldots,h_{n-1})=\nu(\phi_1^{+}, \ldots,\phi_{n-1}^{+}) =\mu_n^{+} >0. \leqno(2.1.8)
$$
Furthermore, if there are only $n-1$ strictly positive eigenvalues, then since $H$ is infinite dimensional there are either infinitely many strictly negative eigenvalues converging to zero, or zero is an eigenvalue with infinite multiplicity. In either of these cases  
we have
$$
\sup_{||x||=1, (x,\phi_j^{+})=0, j=1,\ldots,n-1,(x,\theta)=0, \theta \in \Theta} (Ax,x)=0, \leqno(2.1.9)
$$
where $\Theta$ is any finite set of orthonormal vectors each of which is orthogonal to all $\phi_j^{+}, j=1,\ldots,n-1$. For $\mu_n^{-}$, the result is analogous, but the roles of sup and inf must be interchanged.
\bigskip

An application of the previous theorem provides a comparison of the eigenvalues of several operators. Some results of this type, due to H. Weyl and R. Courant, are referenced on page 238 of [5]. All appeared between 1911 and 1920.
\bigskip

{\bf Theorem 2.} Let $A_1$ and $A_2$ be compact self-adjoint non-zero operators on $H$ and let $A=A_1 + A_2$. Denote the $n^{th}$ non-negative eigenvalues of $A_1,~A_2,$ and $A$ by 
$\mu_{n,A_{1}}^{+},~\mu_{n,A_{2}}^{+}$, and $\mu_{n,A}^{+}$, respectively, and their $n^{th}$ negative eigenvalues by  $\mu_{n,A_{1}}^{-},~\mu_{n,A_{2}}^{-}$, and $\mu_{n,A}^{-}$. Then, for all integers $p,q \ge 1$ 
$$
\mu_{p+q-1,A}^{+} \le \mu_{p,A_{1}}^{+} + \mu_{q,A_{2}}^{+} \leqno (2.1.10)
$$
provided the corresponding strictly positive eigenvalues  in (2.1.10) 
all exist, and 
$$
\mu_{p+q-1,A}^{-} \ge \mu_{p,A_{1}}^{-} + \mu_{q,A_{2}}^{-} \leqno(2.1.11)
$$
provided the corresponding strictly negative eigenvalues in (2.1.11) 
all exist.
In particular, if $q=1$ and the various eigenvalues all exist, then $\max\{ \mu_{1,A_{2}}^{+}, | \mu_{1,A_{2}}^{-}|\} =||A_2||$, 
$$
\mu_{p,A}^{+} \le \mu_{p,A_{1}}^{+} + \mu_{1,A_{2}}^{+} \le \mu_{p,A_{1}}^{+} + ||A_2||, \leqno(2.1.12)
$$
and 
$$
\mu_{p,A}^{-} \ge \mu_{p,A_{1}}^{-} + \mu_{1,A_{2}}^{-}\ge \mu_{p,A_{1}}^{-} - ||A_2||. \leqno(2.1.13)
$$
\bigskip

{\bf Proof.} If the strictly positive eigenvalues appearing in (2.1.10)  
all exist and 
$$
\Gamma_{p,q} =\Phi_{p,A_1} \cup \Phi_{q,A_2},
$$
where
$$
\Phi_{p,A_1} =\{ \phi_{j, A_1}^{+}: 1\le j \le p-1\} ~{\rm{and}}~ \Phi_{q,A_2} =\{\phi_{j, A_2}^{+}:1 \le j \le q-1\},
$$
are the corresponding eigenvectors repeated according to  their multiplicity, then using  Theorem 1 in the first inequality that follows we have
$$
\mu_{p+q-1,A}^{+} \le \sup_{||x||=1, (x,\phi)=0 ~\forall \phi \in \Gamma_{p,q}} (Ax,x) \leqno(2.1.14)
$$
$$
 \le  \sup_{||x||=1, (x,\phi)=0 ~\forall \phi \in \Gamma_{p,q}} (A_1x,x) + \sup_{||x||=1, (x,\phi)=0 ~\forall \phi \in \Gamma_{p,q}} (A_2x,x). 
$$
$$
 \le  \sup_{||x||=1, (x,\phi)=0 ~\forall \phi \in \Phi_{p,A_1}} (A_1x,x) + \sup_{||x||=1, (x,\phi)=0 ~\forall \phi \in \Phi_{q,A_2}} (A_2x,x). 
$$
Hence, by definition of $\mu_{p,A_1}^{+}$ and $\mu_{q,A_2}^{+}$ (2.1.14)  
implies (2.1.10).

If the strictly negative eigenvalues appearing in (2.1.11) 
all exist, then since $\mu_{k,A}^{-}= -\mu_{k,-A}^{+}$ for all $k \ge 1$, it follows from (2.1.10)  
applied to $-A$ that
$$
\mu_{p+q-1,-A}^{+} \le \mu_{p,-A_{1}}^{+} + \mu_{q,-A_{2}}^{+}. \leqno(2.1.15)
$$
Multiplying (2.1.15)
by minus one, then implies (2.1.11). 

Since (2.1.12) 
and (2.1.13) 
are special cases of (2.1.10)
and (2.1.11)  
and that  $\max\{ \mu_{1,A_{2}}^{+}, | \mu_{1,A_{2}}^{-}|\} =||A_2||$, Theorem 2 is proved.
\bigskip

{\bf Remark 1.} Letting $A_1= A- A_2= A + (-A_2)$ and applying (2.1.12) 
with $q=1$, we have
$$
\mu_{p,A_{1}}^{+} \le \mu_{p,A}^{+} + \mu_{1,-A_{2}}^{+}. \leqno(2.1.16)
$$
Combining (2.1.14) 
and (2.1.16)  
now implies
$$
-\mu_{1,-A_{2}}^{+} \le |\mu_{p,A}^{+} -  \mu_{p,A_{1}}^{+}| \le \mu_{1,A_{2}}^{+}, \leqno(2.1.17)
$$
and therefore
$$
-||A_2|| \le |\mu_{p,A}^{+} -  \mu_{p,A_{1}}^{+}| \le||A_2||, \leqno(2.1.18)
$$
since $\max\{\mu_{1,-A_{2}}^{+},\mu_{1,A_{2}}^{+}\} \le||A_2||$.  Again letting $A_1= A- A_2= A + (-A_2)$ and applying (2.1.13) 
with $q=1$, we have
$$
\mu_{p,A_{1}}^{-} \ge \mu_{p,A}^{-} + \mu_{1,-A_{2}}^{-}. \leqno(2.1.19) 
$$
 Combining (2.1.15)  
 and (2.1.19) 
implies
$$
 \mu_{1,A_{2}}^{-} \le \mu_{p,A}^{-} - \mu_{p,A_{1}}^{-} \le - \mu_{1,-A_{2}}^{-}, \leqno(2.1.20)
$$
and since $ -\mu_{1,-A_{2}}^{-}=  \mu_{1,A_{2}}^{+}$ and $ \mu_{1,A_{2}}^{-} =- \mu_{1,-A_{2}}^{+} $ it follows that
$$
| \mu_{p,A}^{-} - \mu_{p,A_{1}}^{-} | \le  \max\{\mu_{1,-A_{2}}^{+},\mu_{1,A_{2}}^{+}\} \le||A_2||. \leqno(2.1.21)
$$
\bigskip

{\bf 2.2. Eigenvalue Comparison for Self-Adjoint Compact Operators.} Here we use some of the ideas in Section 2.1, but our approach is more direct, and we also deal with the Hausdorff distance between the spectrums

As before we assume $A,A_1,A_2$ are compact self-adjoint operators on the infinite dimensional Hilbert space $H$ and that $A=A_1+A_2$. Since 
$$
\sup_{x\in H,||x||=1}|(A_2x,x)| =||A_2||,
$$
then for every $x \in H$
$$
(Ax,x) \le (A_1x,x) + ||A_2|| ~{\rm {and}}~ (Ax,x) \ge (A_1x,x) - ||A_2||. \leqno(2.2.1)
$$
Hence if $\mu_{p,A}^{+}$ and $\mu_{p,A_1}^{+}$ both exist for some $p \ge 1$, and
 $$
 \Phi_{p, A_1} =\{\phi_{j,A_1}^{+}: 1 \le j \le p-1\}
 $$
 are the eigenvectors corresponding to the strictly positive eigenvalues $\{\mu_{j,A_1}^{+}: 1 \le j \le p-1 \}$ of $A_1$, then by definition
 $$
\mu_{p,A_1}^{+}=  \sup_{||x||=1, (x,\phi)=0~\forall \phi \in \Phi_{p,A_1}} (A_1x,x),
$$ 
and by Theorem 1
$$
\sup_{||x||=1, (x,\phi)=0~\forall~ \phi \in \Phi_{p,A_1}} (Ax,x) \ge \mu_{p,A}^{+}.
$$ 
Therefore,
$$
 \mu_{p,A}^{+} \le  \mu_{p,A_1}^{+} +||A_2||. \leqno(2.2.2)
$$
Similarly, if $\Phi_{p,A}= \{\phi_{j,A}^{+}:1\le j \le p-1\}$ are the eigenvectors corresponding to the strictly positive eigenvalues $\{\mu_{j,A}^{+}:1 \le j \le p-1\}$, then
by definition
 $$
\mu_{p,A}^{+}=  \sup_{||x||=1, (x,\phi)=0~\forall  \phi \in \Phi_{p,A}} (Ax,x),
$$ 
and by Theorem 1
$$
\sup_{||x||=1, (x,\phi)=0~\forall \phi \in \Phi_{p,A}} (A_1x,x) \ge \mu_{p,A_1}^{+}.
$$ 
Therefore,
$$
 \mu_{p,A}^{+} \ge  \mu_{p,A_1}^{+} -||A_2||, \leqno(2.2.3)
$$
and combining (2.2.2)  
and (2.2.3)  
$$
-||A_2|| \le  \mu_{p,A}^{+} - \mu_{p,A_1}^{+} \le ||A_2||, \leqno(2.2.4)
$$
which is a more direct way to obtain (2.1.18). 

For $B$ a compact self-adjoint operator on $H$, we let $S(B)$ denote the eigenvalues of $B$,  $S(B)^{+}$ the strictly positive eigenvalues of $B$ and $S(B)^{-}$ the strictly negative eigenvalues of $B$. Then by Theorem 1 it can easily be seen that zero may or may not be an eigenvalue of $B$, but if $0 \notin S(B)$, it is always a limit point of either $S(B)^{+}$ or $S(B)^{-}$.
\bigskip

{\bf Lemma 1.} Let  $A,A_1,A_2$ be compact self-adjoint operators on the infinite dimensional Hilbert space $H$ such that $A=A_1+A_2.$ If $\delta>||A_2||$, then
$$
S(A)^{+} \subseteq (S(A_1)^{+} \cup \{0\}) + (-\delta,\delta), \leqno(2.2.5)
$$
and
$$ 
S(A_1)^{+} \subseteq (S(A)^{+} \cup \{0\}) + (-\delta,\delta), \leqno(2.2.6)
$$
which together imply the Hausdorff distance between the sets $S(A)^{+} \cup\{0\})$ and $S(A_1)^{+} \cup\{0\})$ is less than or equal to $||A_2||$.

\bigskip

{\bf Proof.} Let $A_1$ have $r$ strictly positive eigenvalues. Then $r=0$,  $1 \le r \le \infty$ and $r=\infty$, and we also observe that when
$r< \infty$ , then $H$ infinite dimensional implies $A_1$ has infinitely many strictly negative eigenvalues or zero is an eigenvalue corresponding to infinitely many orthonormal eigenvectors. 

If $r=0$, then by (2.1.8)  
and (2.1.9)  
$$
\sup_{||x||=1}(A_1x,x)=0,
$$
and (2.2.1)
implies
$$
\sup_{||x||=1}(Ax,x) \le  0+ ||A_2|| \leqno(2.2.7)
$$
and
$$
\sup_{||x||=1}(Ax,x) \ge  0- ||A_2||, \leqno(2.2.8)
$$
which implies every strictly positive eigenvalue of $A$ (if any exist) must lie in the open interval $(-\delta,\delta)$ whenever $\delta>||A_2||$. Thus, (2.2.5)
hold if $r=0$.
\bigskip

If $r =\infty$ and $1 \le p <\infty$, or $1\le r <\infty$ and $1\le p \le r$, then for all such $r$ and $p$ whenever the strictly positive eigenvalue $\mu_{p,A}^{+}$ exists we have from (2.2.4) 
that
$$
\mu_{p,A}^{+} \in \mu_{p,A_{1}}^{+} +(-\delta,\delta). \leqno(2.2.9)
$$
\bigskip

If $1\le r <\infty$  and $p > r$, then for every strictly positive eigenvalue $\mu_{p,A}^{+}$ that exists, we next show that
$$
0 \le \mu_{p,A}^{+} \le ||A_2||, \leqno(2.2.10)
$$
which implies
$$
\mu_{p,A}^{+} \in (-\delta,\delta). \leqno(2.2.11)
$$
Combining (2.2.9), 
(2.2.11),  
and the case $r=0$, we again have (2.2.5) 
once (2.2.10) 
is verified. 
\bigskip

To verify (2.2.10)  
we recall $\Phi_{r,A_{1}}$ as above and define
$$
\Phi_{p}= \Phi_{r,A_{1}} \cup \Theta_{r,p}
$$
where $\Theta_{r,p}$ is a set of $p-r$ orthonormal vectors which are orthogonal to all the eigenvectors in $\Phi_{r,A_{1}}$. Then,
$$
\mu_{p,A}^{+} \le \sup_{||x||=1, (x,v)=0 ~\forall v \in \Phi_p} (Ax,x) \le  \sup_{||x||=1, (x,v)=0 ~\forall v \in \Phi_p} (A_1x,x)+ ||A_2|| \le 0 +||A_2|| \leqno(2.2.12)
$$
where the first inequality in (2.2.12) 
follows from Theorem 1 and the second from (2.1.9). 

Hence (2.2.5)  
is proven, and to verify (2.2.6) 
we note that $A_1=A +(-A_2)$ and then simply repeat the previous argument starting with (2.2.1), 
interchanging $A$ and $A_1$ and replacing $A_2$ with $-A_2$. Thus the lemma is proven.
\bigskip

The analogue of Lemma 1 for the strictly negative eigenvalues for these operators as follows.
\bigskip

{\bf Lemma 2.} Let  $A,A_1,A_2$ be compact self-adjoint operators on the infinite dimensional Hilbert space $H$ such that $A=A_1+A_2.$ If $\delta>||A_2||$, then
$$
S(A)^{-} \subseteq (S(A_1)^{-} \cup \{0\}) + (-\delta,\delta), \leqno(2.2.13)
$$
and
$$ 
S(A_1)^{-} \subseteq (S(A)^{-} \cup \{0\}) + (-\delta,\delta), \leqno(2.2.14)
$$
which together imply the Hausdorff distance between the sets $S(A)^{-} \cup\{0\})$ and $S(A_1)^{-} \cup\{0\})$ is less than or equal to $||A_2||$.
\bigskip

{\bf Proof.} Since the strictly negative eigenvalues of $A$ and  $A_1$ are the strictly positive eigenvalues of $-A$ and $A_1$ multiplied by $-1$, the proof of Lemma 1 applied to
$$
-A= -A_1 + (-A_2)
$$
establishes (2.2.13)  
and (2.2.14).  
Hence the lemma is proven.

Combining Lemma 1 and Lemma 2, and keeping in mind that $H$ infinite dimensional implies zero is a limit point of the eigenvalues (or actually an eigenvalue) of both $A$ and $A_1$, we have the following theorem.
\bigskip

{\bf Theorem 3.}  Let  $A,A_1,A_2$ be compact self-adjoint operators on the infinite dimensional Hilbert space $H$ such that $A=A_1+A_2.$ If $\delta>||A_2||$, then
$$
S(A) \subseteq S(A_1) + (-\delta,\delta)~{\rm{and}}~S(A_1) \subseteq S(A) + (-\delta,\delta). \leqno(2.2.15)
$$
which together imply the Hausdorff distance between the sets $S(A)$ and $S(A_1)$ is less than or equal to $||A_2||$.
\bigskip

{\bf References.}
\bigskip

[1] J.N. Franklin, Matrix Theory, Prentice-Hall, Englewood Cliffs, N.J., 1968. 

[2] L. Gross, Potential Theory on Hilbert Spaces, 1967, J. Functional Analysis, Vol. 1, 123-181.

[3] J. Kuelbs and J. Zinn, Limits for Partial Maxima of Gaussian Random Vectors, submitted for publication in the Journal of Theoretical Probability.

[4] Yu. V. Prokhorov, Convergence of Random Processes and Limit Theorem in Probability, Theor. Probab. Appl. 1956, Vol. 1, 157-214.

[5]  F. Riesz, B. Sz.-Nagy, Functional Analysis, Frederick Ungar Publishing, New York, 1955.
\bigskip

\end{document}